\documentclass[11pt]{article} 
\usepackage{amsmath,amssymb,amsthm}
\usepackage{fancybox}  
\usepackage{graphicx}
\usepackage{color}
\usepackage{colortab}
\usepackage{colortbl}
\usepackage{mathdots}

\allowdisplaybreaks

\begin{document}

\def\fl#1{\left\lfloor#1\right\rfloor}
\def\cl#1{\left\lceil#1\right\rceil}
\def\ang#1{\left\langle#1\right\rangle}
\def\stf#1#2{\left[#1\atop#2\right]} 
\def\sts#1#2{\left\{#1\atop#2\right\}}
\def\eul#1#2{\left\langle#1\atop#2\right\rangle}
\def\N{\mathbb N}
\def\Z{\mathbb Z}
\def\R{\mathbb R}
\def\C{\mathbb C}
\newcommand{\ctext}[1]{\raise0.2ex\hbox{\textcircled{\scriptsize{#1}}}}

\newtheorem{theorem}{Theorem}
\newtheorem{Prop}{Proposition}
\newtheorem{Cor}{Corollary}
\newtheorem{Lem}{Lemma}
\newtheorem{Sublem}{Sublemma}
\newtheorem{Def}{Definition}  
\newtheorem{Conj}{Conjecture}

\newenvironment{Rem}{\begin{trivlist} \item[\hskip \labelsep{\it
Remark.}]\setlength{\parindent}{0pt}}{\end{trivlist}}

\title{$p$-numerical semigroup of generalized Fibonacci triples}

\author{
Takao Komatsu 
\\
\small Department of Mathematical Sciences, School of Science\\[-0.8ex]
\small Zhejiang Sci-Tech University\\[-0.8ex]
\small Hangzhou 310018 China\\[-0.8ex]
\small \texttt{komatsu@zstu.edu.cn}\\\\
Shanta Laishram\\
\small Department of Mathematics \\[-0.8ex]
\small Theoretical Statistical \& Mathematics Unit (SMU)\\[-0.8ex]
\small Indian Statistical Institute\\[-0.8ex]
\small New Delhi 110016, India\\[-0.8ex]
\small \texttt{shanta@isid.ac.in}\\\\
Pooja Punyani\\
\small Department of Applied Sciences\\[-0.8ex] 
\small The NorthCap University\\[-0.8ex]  
\small Gurugram 122001 Haryana, India\\[-0.8ex] 
\small \texttt{poojapunyani@ncuindia.edu}  
}

\date{
\small MR Subject Classifications: Primary 11D07; Secondary 20M14, 05A17, 05A19, 11D04, 11B68, 11P81 
}

\maketitle
 
\begin{abstract} 
For a nonnegative integer $p$, we give explicit formulas for the $p$-Frobenius number and the $p$-genus of generalized Fibonacci numerical semigroups. Here, the $p$-numerical semigroup $S_p$ is defined as the set of integers whose nonnegative integral linear combinations of given positive integers $a_1,a_2,\dots,a_k$ are expressed more than $p$ ways. When $p=0$, $S_0$ with the $0$-Frobenius number and the $0$-genus is the original numerical semigroup with the Frobenius number and the genus. In this paper, we consider the $p$-numerical semigroup involving Jacobsthal polynomials, which include Fibonacci numbers as special cases. We can also treat with the Jacobsthal-Lucas polynomials, including Lucas numbers accordingly. One of the applications on the $p$-Hilbert series is mentioned. 
\end{abstract}

\section{Introduction}  

Given the set of positive integers $A:=\{a_1,a_2,\dots,a_k\}$ ($k\ge 2$), for a nonnegative integer $p$, let $S_p$ be the set of integers whose nonnegative integral linear combinations of given positive integers $a_1,a_2,\dots,a_k$ are expressed in more than $p$ ways.  
For some backgrounds of the number of representations, see, e.g., \cite{Binner20,Cayley,Ko03,sy1857,tr00}.  
For a set of nonnegative integers $\mathbb N_0$, the set $\mathbb N_0\backslash S_p$ is finite if and only if $\gcd(a_1,a_2,\dots,a_k)=1$.  Then there exists the largest integer $g_p(A):=g(S_p)$ in $\mathbb N_0\backslash S_p$, which is called the {\it $p$-Frobenius number}. The cardinality of $\mathbb N_0\backslash S_p$ is called the {\it $p$-genus} and is denoted by $n_p(A):=n(S_p)$.   
The sum of the elements in $\mathbb N_0\backslash S_p$ is called the {\it $p$-Sylvester sum} and is denoted by $s_p(A):=s(S_p)$.  
 This kind of concept is a generalization of the famous Diophantine problem of Frobenius (\cite{RA05,sy1882,sy1884}) since $p=0$ is the case when the original Frobenius number $g(A)=g_0(A)$, the genus $n(A)=n_0(A)$   
and the Sylvester sum $s(A)=s_0(A)$   
are recovered. We can call $S_p$ the {\it $p$-numerical semigroup}. Strictly speaking, when $p\ge 1$, $S_p$ does not include $0$ since the integer $0$ has only one representation, so it satisfies simply additivity and the set $S_p\cup\{0\}$ becomes a numerical semigroup. For numerical semigroups, we refer to \cite{ADG,RG1,RG2}.  

There exist different extensions of the Frobenius number and genus, even in terms of the number of representations called denumerant. For example, some considered $S_p^\ast$ as the set of integers whose nonnegative integral linear combinations of given positive integers $a_1,a_2,\dots,a_k$ are expressed in {\it exactly} $p$ ways (see, e.g., \cite{BDFHKMRSS,FS11}). Then, 
the corresponding $p$-Frobenius number $g_p^\ast(A)$ is the largest integer that has {\it\underline{exactly $p$ distinct}} representations. However, in this case, $g_p^\ast(A)$ is not necessarily increasing as $p$ increases. For example, when $A:=\{2,5,7\}$, $g_{17}^\ast(2,5,7)=43>g_{18}^\ast(2,5,7)=42$. In addition, for some $p$, $g_p^\ast$ may not exist. For example, $g_{22}^\ast(2,5,7)$ does not exist because there is no positive integer whose number of representations is exactly $22$.  
The $p$-genus may be also defined in different ways. For example, $n_p^\ast(A)$ can be defined the cardinality of $[\ell_p(A),g_p(A)+1]\backslash S_p(A)$, where $\ell_p(A)$ is the least element of $S_p(A)$. However, in our definition of $n_p(A)$ as the cardinality of   $[0,g_p(A)+1]\backslash S_p(A)$, we can use the convenient formula arising from the $p$-Ap\'ery set in order to obtain $n_p(A)$. See the next section.  
      
In \cite{Matt09}, numerical semigroups generated by $\{a,a+b,a F_{k-1}+b F_k\}$ are considered. Using a technique of Johnson \cite{Johnson60}, the Frobenius numbers of such semigroups are found as a generalization of the result by Marin et. al. \cite{MAR}. 

In this paper, for a positive integer $v$, we treat with Jacobsthal polynomials $J_n(v)$, defined by the recurrence relation $J_n(v)=J_{n-1}(v)+v J_{n-2}(v)$ ($n\ge 2$) with $J_0(v)=0$ and $J_1(v)=1$ (see, e.g., \cite[Chapter 44]{Koshy2}). When $v=1$, $F_n=J_n(1)$ are Fibonacci numbers. When $v=2$, $J_n=J_n(2)$ are Jacobsthal numbers. Then, 
we give explicit formulas of $p$-Frobenius numbers for $A:=\{a,v a+b,v a J_{k-1}(v)+b J_k(v)\}$, where $a$ and $b$ are positive integers with $\gcd(a,b)=1$ and $a,k\ge 3$. If $a=J_i(v)$ and $b=J_{i+1}(v)$, then by $v J_i(v)J_{k-1}(v)+J_{i+1}(v) J_k(v)=J_{i+k}(v)$, we get $A=\{J_i(v),J_{i+2}(v),J_{i+k}(v)\}$. So, the results in \cite{KP} are recovered as special cases. In addition, if $v=1$, the results in \cite{MAR,KY} are recovered as special cases. 
  
For $k=2$, that is, the case of two variable, closed formulas are explicitly given for $g_0(A)$ (\cite{sy1884}), $n_0(A)$ (\cite{sy1882}) and $s_0(A)$ (\cite{bs93}; its extension \cite{ro94}). 
However, for $k\ge 3$, the Frobenius number cannot be given by any set of closed formulas which can be reduced to a finite set of certain polynomials (\cite{cu90}). For $k=3$, various algorithms have been devised for finding the Frobenius number (\cite{Fel06,Johnson60,RGS04}). Some inexplicit formulas for the Frobenius number in three variables can be seen in \cite{tr17}.   
The need to care enough is that even in the original case of $p=0$, it is very difficult to give a closed explicit formula of any general sequence for three or more variables (see, e.g., \cite{RR18,RGS04,RBT2015,RBT2017}). Indeed, it is even more difficult when $p>0$. However, finally, we have succeeded in obtaining $p$-Frobenius number in triangular numbers \cite{Ko22a} and repunits \cite{Ko22b} as well as Fibonacci and Lucas triplets \cite{KY} and Jacobsthal triples \cite{KP} quite recently.

\section{Preliminaries}  

We introduce the Ap\'ery set (see \cite{Apery}) below in order to obtain the formulas for $g_p(A)$, $n_p(A)$ and $s_p(A)$ technically. Without loss of generality, we assume that $a_1=\min(A)$. 

\begin{Def}  
Let $p$ be a nonnegative integer. For a set of positive integers $A=\{a_1,a_2,\dots,a_k\}$ with $\gcd(A)=1$ and $a_1=\min(A)$ we denote by 
$$
{\rm Ap}_p(A)={\rm Ap}_p(a_1,a_2,\dots,a_k)=\{m_0^{(p)},m_1^{(p)},\dots,m_{a_1-1}^{(p)}\}\,, 
$$  
the $p$-Ap\'ery set of $A$, where $m_i^{(p)}$ is the least positive integer of $S_p(A)$, satisfying $m_i^{(p)}\equiv i\pmod{a_1}$ $(0\le i\le a_1-1)$ (that is, $m_i^{(p)}\in S_p(A)$ and $m_i^{(p)}-a_1\not\in S_p(A)$). 
Note that $m_0^{(0)}$ is defined to be $0$.  
\label{apery} 
\end{Def}  

\noindent 
It follows that for each $p$, 
$$
{\rm Ap}_p(A)\equiv\{0,1,\dots,a_1-1\}\pmod{a_1}\,. 
$$  

Even though it is hard to find any explicit form of $g_p(A)$ as well as $n_p(A)$ and $s_p(A)$ $k\ge 3$, by using convenient formulas established in \cite{Ko22,Ko-p}, we can obtain such values for some special sequences $(a_1,a_2,\dots,a_k)$ after finding any regular structure of $m_j^{(p)}$ is enough hard in general. One convenient formula is on the power sum 
$$
s_p^{(\mu)}(A):=\sum_{n\in\mathbb N_0\backslash S_p(A)}n^\mu
$$  
by using Bernoulli numbers $B_n$ defined by the generating function 
$$
\frac{x}{e^x-1}=\sum_{n=0}^\infty B_n\frac{x^n}{n!}\,, 
$$
and another convenient formula is on the weighted power sum (\cite{KZ0,KZ}) 
$$
s_{\lambda,p}^{(\mu)}(A):=\sum_{n\in \mathbb N_0\backslash S_p(A)}\lambda^n n^\mu 
$$  
by using Eulerian numbers $\eul{n}{m}$ appearing in the generating function 
$$ 
\sum_{k=0}^\infty k^n x^k=\frac{1}{(1-x)^{n+1}}\sum_{m=0}^{n-1}\eul{n}{m}x^{m+1}\quad(n\ge 1)
$$ 
with $0^0=1$ and $\eul{0}{0}=1$. Here, $\mu$ is a nonnegative integer and $\lambda\ne 1$. From these convenient formulas, many useful expressions are yielded as special cases. 
Some useful ones are given as follows.  The formulas (\ref{mp-n}) and (\ref{mp-s}) are entailed from $s_{\lambda,p}^{(0)}(A)$ and $s_{\lambda,p}^{(1)}(A)$, respectively.  

\begin{Lem}  
Let $k$, $p$ and $\mu$ be integers with $k\ge 2$, $p\ge 0$ and $\mu\ge 1$.  
Assume that $\gcd(a_1,a_2,\dots,a_k)=1$.  We have 
\begin{align}  
g_p(a_1,a_2,\dots,a_k)&=\left(\max_{0\le j\le a_1-1}m_j^{(p)}\right)-a_1\,, 
\label{mp-g}
\\  
n_p(a_1,a_2,\dots,a_k)&=\frac{1}{a_1}\sum_{j=0}^{a_1-1}m_j^{(p)}-\frac{a_1-1}{2}\,,  
\label{mp-n}
\\
s_p(a_1,a_2,\dots,a_k)&=\frac{1}{2 a_1}\sum_{j=0}^{a_1-1}\bigl(m_j^{(p)}\bigr)^2-\frac{1}{2}\sum_{j=0}^{a_1-1}m_j^{(p)}+\frac{a_1^2-1}{12}\,.
\label{mp-s}
\end{align} 
\label{lem-mp}
\end{Lem} 

\noindent 
{\it Remark.}  
When $p=0$, the formulas (\ref{mp-g}), (\ref{mp-n}) and (\ref{mp-s}) reduce to the formulas by  
Brauer and Shockley \cite[Lemma 3]{bs62}, Selmer \cite[Theorem]{se77}, and Tripathi \cite[Lemma 1]{tr08}\footnote{There was a typo, but it was corrected in \cite{PT}.}, respectively: 
\begin{align*}   
g(a_1,a_2,\dots,a_k)&=\left(\max_{0\le j\le a_1-1}m_j\right)-a_1\,,\\  
n(a_1,a_2,\dots,a_k)&=\frac{1}{a_1}\sum_{j=0}^{a_1-1}m_j-\frac{a_1-1}{2}\,,\\  
s(a_1,a_2,\dots,a_k)&=\frac{1}{2 a_1}\sum_{j=0}^{a_1-1}(m_j)^2-\frac{1}{2}\sum_{j=0}^{a_1-1}m_j+\frac{a_1^2-1}{12}\,,   
\end{align*} 
where $m_j=m_j^{(0)}$ ($1\le j\le a_1-1$) with $m_0=m_0^{(0)}=0$.

\section{Main results}  

Determine integers $q$ and $r$ by $a=q J_k(v)+r$ with $0\le r<J_k(v)$.  
The function $\fl{x}$ denotes the largest integer that does not exceed $x$. 

\begin{theorem}
Let $a$ and $b$ be positive integers with $a\ge 3$ and $\gcd(a,b)=1$. Then, 
for a positive integer $k\ge 3$, and $0\le p\le\fl{a/J_k(v)}$ we have 
\begin{align*}
&g_p\bigl(a,v a+b,v a J_{k-1}(v)+b J_k(v)\bigr)\\
&=\begin{cases}
(a-1)b+a\bigl(v(r-1)-1\bigr)+\frac{v a(a-r)J_{k-1}(v)}{J_k(v)}+p\bigl(v a J_{k-1}(v)+b J_k(v)\bigr)&\\
\qquad\qquad\qquad\text{if $a<J_k(v)$ or $(v a+b)r>v a\bigl(J_k(v)-J_{k-1}(v)\bigr)$};&\\ 
(a-r-1)b+v a(J_k(v)-J_{k-1}(v)-1)-a+\frac{v a(a-r)J_{k-1}(v)}{J_k(v)}&\\
+p\bigl(v a J_{k-1}(v)+b J_k(v)\bigr)\quad\text{if $a\ge J_k(v)$ and $(v a+b)r<v a\bigl(J_k(v)-J_{k-1}(v)\bigr)$}\,,& 
\end{cases}
\end{align*} 
where $r=a-\fl{a/J_k(v)}J_k(v)$. 
\label{th1}
\end{theorem}  

For example, if $k=3$ and $v=1$, then for $0\le p\le\fl{a/2}$ we have 
\begin{align*}
&g_p(a,a+b,a F_{2}+b F_3)\\
&=\begin{cases}
(a-1)b+\dfrac{a(a-3)}{2}+p(a+2 b)&\text{if $a$ is odd};\\ 
(a-1)b+\dfrac{a(a-2)}{2}+p(a+2 b)&\text{if $a$ is even}\,. 
\end{cases}
\end{align*}

\subsection{The case $p=0$} 

In this triple $\{a,v a+b,v a J_{k-1}(v)+b J_k(v)\}$, we can use the similar framework to that in \cite{KP} to construct the elements of the $p$-Ap\'ery set. It is very important to see that such a framework is not always possible. See \cite{Ko22a} for a different structure. No structure has been analyzed for most other triples, so no explicit formula has been found for them.    

Consider the expression 
$$
t_{y,z}:=y(v a+b)+z\bigl(v a J_{k-1}(v)+b J_k(v)\bigr)\,. 
$$
We see that $q=\fl{a/J_k(v)}$. 
Then all the elements in $0$-Ap\'ery set are represented as in Table \ref{tb:ap0}. 

\begin{table}[htbp]
  \centering
\begin{tabular}{ccccc}
\cline{1-2}\cline{3-4}\cline{5-5}
\multicolumn{1}{|c}{$t_{0,0}$}&$\cdots$&$\cdots$&$\cdots$&\multicolumn{1}{c|}{$t_{J_k(v)-1,0}$}\\
\multicolumn{1}{|c}{$t_{0,1}$}&$\cdots$&$\cdots$&$\cdots$&\multicolumn{1}{c|}{$t_{J_k(v)-1,1}$}\\
\multicolumn{1}{|c}{$\vdots$}&&&&\multicolumn{1}{c|}{$\vdots$}\\
\multicolumn{1}{|c}{$t_{0,q-1}$}&$\cdots$&$\cdots$&$\cdots$&\multicolumn{1}{c|}{$t_{J_k(v)-1,q-1}$}\\
\cline{4-5}
\multicolumn{1}{|c}{$t_{0,q}$}&$\cdots$&\multicolumn{1}{c|}{$t_{r-1,q}$}&&\\
\cline{1-2}\cline{2-3}
\end{tabular}
  \caption{${\rm Ap}_0(a,v a+b,v a J_{k-1}(v)+b J_k(v))$}
  \label{tb:ap0}
\end{table}

Since $t_{i+1,j}-t_{i,j}\equiv b\pmod a$ and $t_{0,j+1}-t_{J_k(v)-1,j}\equiv b\pmod a$, the sequence 
\begin{align*} 
&t_{0,0},t_{1,0},\dots,t_{J_k(v)-1,0},t_{0,1},t_{1,1},\dots,t_{J_k(v)-1,1},\dots\\
&\quad t_{0,q-1},t_{1,q-1},\dots,t_{J_k(v)-1,q-1},t_{0,q},\dots,t_{r-1,q}
\end{align*} 
is equivalent to the sequence $\{\ell b\pmod a\}_{\ell=0}^{a-1}$ in a way that keeps this order completely. 
Since $\gcd(a,b)=1$ (otherwise, $\gcd(A)>1$), all the elements appeared in Table \ref{tb:ap0} constitute a complete residue system modulo $a\bigl(=q J_k(v)+r\bigr)$. 

It is clear that the largest element in ${\rm Ap}_0(A)$, where $A:=\{a,v a+b,v a J_{k-1}(v)+b J_k(v)\}$, is $t_{r-1,q}$ or $t_{J_k(v)-1,q-1}$. 
If $a<J_k(v)$, by $q=0$, the largest element is $t_{r-1,q}=t_{a-1,q}$. 
Otherwise, by $q>0$, compare two values. The inequality $t_{r-1,q}>t_{J_k(v)-1,q-1}$ holds if and only if  
$(v a+b)r>v a\bigl(J_k(v)-J_{k-1}(v)\bigr)$. 
Hence, if $(v a+b)r>v a\bigl(J_k(v)-J_{k-1}(v)\bigr)$, 
then by Lemma \ref{lem-mp} (\ref{mp-g}) we have 
\begin{align*} 
g_0(A)&=t_{r-1,q}-a\\
&=(a-1)b+a\bigl(v(r-1)-1\bigr)+\frac{v a(a-r)J_{k-1}(v)}{J_k(v)}\,. 
\end{align*} 
If $(v a+b)r<v a\bigl(J_k(v)-J_{k-1}(v)\bigr)$, 
then we have 
\begin{align*} 
g_0(A)&=t_{J_k(v)-1,q-1}-a\\
&=(a-r-1)b+v a(J_k(v)-J_{k-1}(v)-1)-a+\frac{v a(a-r)J_{k-1}(v)}{J_k(v)}\,. 
\end{align*}      
Note that $(v a+b)r\ne v a\bigl(J_k(v)-J_{k-1}(v)\bigr)$ because $\gcd(a,b)=1$.

\subsection{The case $p=1$}  

We assume that $a\ge J_k(v)$ from now on. If $a<J_k(v)$, the situation becomes more and more complicated by requiring a lot of case-by-case discussion for $p\ge 1$. So, the discussion that follows does not apply. 

All the elements in ${\rm Ap}_1(A)$ can be determined from those in ${\rm Ap}_0(A)$. 
Only those elements that have the same residue modulo $a$ as those in the top row of ${\rm Ap}_0(A)$ are arranged in order in the form of filling gaps under the same block. Elements that have the same residue modulo $a$ as the other elements of ${\rm Ap}_0(A)$ are arranged in a row shift up to the immediately adjacent block. It is shown in Table \ref{tb:ap0-1}.

\begin{table}[htbp]
  \centering
\scalebox{0.7}{                  
  \begin{tabular}{cccccccccccc}
\cline{1-2}\cline{3-4}\cline{5-6}\cline{7-8}\cline{9-10}\cline{11-12}
\multicolumn{1}{|c}{\shadowbox{$t_{0,0}$}}&\shadowbox{$t_{1,0}$}&$\cdots$&&$\cdots$&\multicolumn{1}{c|}{\shadowbox{$t_{J_k(v)-1,0}$}}&$t_{J_k(v),0}$&$t_{J_k(v)+1,0}$&$\cdots$&&$\cdots$&\multicolumn{1}{c|}{$t_{2 J_k(v)-1,0}$}\\
\multicolumn{1}{|c}{$t_{0,1}$}&$t_{1,1}$&$\cdots$&&$\cdots$&\multicolumn{1}{c|}{$t_{J_k(v)-1,1}$}&$t_{J_k(v),1}$&$t_{J_k(v)+1,1}$&$\cdots$&&$\cdots$&\multicolumn{1}{c|}{$t_{2 J_k(v)-1,1}$}\\
\multicolumn{1}{|c}{$t_{0,2}$}&$t_{1,2}$&$\cdots$&&$\cdots$&\multicolumn{1}{c|}{$t_{J_k(v)-1,2}$}&$\vdots$&$\vdots$&&&&\multicolumn{1}{c|}{$\vdots$}\\
\multicolumn{1}{|c}{$\vdots$}&$\vdots$&&&&\multicolumn{1}{c|}{$\vdots$}&$t_{J_k(v),q-2}$&$t_{J_k(v)+1,q-2}$&$\cdots$&&$\cdots$&\multicolumn{1}{c|}{\Ovalbox{$t_{2 J_k(v)-1,q-2}$}}\\
\cline{10-11}\cline{11-12}
\multicolumn{1}{|c}{$t_{0,q-1}$}&$t_{1,q-1}$&$\cdots$&&$\cdots$&\multicolumn{1}{c|}{$t_{J_k(v)-1,q-1}$}&$t_{J_k(v),q-1}$&$\cdots$&\multicolumn{1}{c|}{\Ovalbox{$t_{J_k(v)+r-1,q-1}$}}&&&\\
\cline{4-5}\cline{5-6}\cline{7-8}\cline{8-9}
\multicolumn{1}{|c}{$t_{0,q}$}&$\cdots$&\multicolumn{1}{c|}{$t_{r-1,q}$}&$t_{r,q}$&$\cdots$&\multicolumn{1}{c|}{\Ovalbox{$t_{J_k(v)-1,q}$}}&&&&&&\\
\cline{1-2}\cline{3-4}\cline{5-6}
\multicolumn{1}{|c}{$t_{0,q+1}$}&$\cdots$&\multicolumn{1}{c|}{\Ovalbox{$t_{r-1,q+1}$}}&&&&&&&&&\\
\cline{1-2}\cline{2-3}
  \end{tabular}
} 
  \caption{${\rm Ap}_1(A)$ from ${\rm Ap}_0(A)$}
  \label{tb:ap0-1}
\end{table}

This fact is supported by the congruence relations  
\begin{align*}  
t_{y,z}&\equiv t_{y+J_k(v),z-1}\pmod{a}\\
&(0\le y\le J_{k}(v)-1,0\le z\le q-1; 0\le y\le r-1,z=q)\,,\\
t_{y,0}&\equiv t_{y+r,q}\pmod{a}\quad(0\le y\le J_k(v)-r-1)\,,\\
t_{J_k(v)-r+y,0}&\equiv t_{y,q+1}\pmod{a}\quad(0\le y\le r-1)\,. 
\end{align*}
In addition, each element of ${\rm Ap}_1(A)$ has two representations in terms of $a$, $v a+b$ and $v a J_{k-1}(v)+b J_k(v)$, because  
\begin{align*}  
t_{y+J_k(v),z-1}&=\bigl(y+J_k(v)\bigr)(v a+b)+(z-1)\bigl(v a J_{k-1}(v)+b J_k(v)\bigr)\\
&=v\bigl(J_{k}(v)-J_{k-1}(v)\bigr)a+y(v a+b)+z\bigl(v a J_{k-1}(v)+b J_k(v)\bigr)\,,\\
t_{y+r,q}&=(y+r)(v a+b)+q\bigl(v a J_{k-1}(v)+b J_k(v)\bigr)\\
&=\bigl(v a+b-q\bigl(J_{k}(v)-J_{k-1}(v)\bigr)\bigr)a+y(v a+b)\,,\\
t_{y,q+1}&=y(v a+b)+(q+1)\bigl(v a J_{k-1}(v)+b J_k(v)\bigr)\\
&=\bigl(v a+b-(q+1)v^2 J_{k-2}(v)\bigr)a+\bigl(y+J_k(v)-r\bigr)(v a+b)\,. 
\end{align*}
Notice that $(v a+b-q\bigl(J_{k}(v)-J_{k-1}(v)\bigr)>0$ and $v a+b-(q+1)v^2 J_{k-2}(v)>0$ because $a=q J_k(v)+r$ and $J_k(v)=J_{k-1}(v)+v J_{k-2}(v)$.   

There are four candidates to take the largest value in ${\rm Ap}_1(A)$: 
$$
t_{r-1,q+1},\quad t_{J_k(v)-1,q},\quad t_{J_k(v)+r-1,q-1},\quad t_{2 J_k(v)-1,q-2}\,. 
$$  
However, since $2 v a J_{k-1}(v)+b J_k(v)>v a J_k(v)$, we can 
see that $t_{r-1,q+1}>t_{J_k(v)+r-1,q-1}$ and $t_{J_k(v)-1,q}>t_{2 J_k(v)-1,q-2}$. In addition, $t_{r-1,q+1}>t_{J_k(v)-1,q}$ if and only if   
$(v a+b)r>v a\bigl(J_k(v)-J_{k-2}(v)\bigr)$.   
Hence, if $(v a+b)r>a\bigl(J_{k}(v)-J_{k-1}(v)\bigr)$, then by Lemma \ref{lem-mp} (\ref{mp-g}) we have  
\begin{align*} 
g_1(A)&=t_{r-1,q+1}-a\\
&=(a-1)b+a\bigl(v(r-1)-1\bigr)+\dfrac{v a(a-r)J_{k-1}(v)}{J_k(v)}+\bigl(v a J_{k-1}(v)+b J_k(v)\bigr)\,. 
\end{align*} 
If $(v a+b)r<a\bigl(J_{k}(v)-J_{k-1}(v)\bigr)$, then we have 
\begin{align*} 
g_1(A)&=t_{J_k(v)-1,q}-a\\
&=(a-r-1)b+v a(J_k(v)-J_{k-1}(v)-1)-a\\
&\quad +\frac{v a(a-r)J_{k-1}(v)}{J_k(v)}+(v a J_{k-1}(v)+b J_k(v))\,. 
\end{align*}

\subsection{The case $p\ge 2$} 

When $p\ge 2$, in a similar manner, each element of ${\rm Ap}_p(A)$ is determined by the corresponding element with the same residue modulo $a$ in ${\rm Ap}_{p-1}(A)$. 
In each block with a lateral length of $J_k(v)$, the elements in the top row in ${\rm Ap}_{p-1}(A)$ are arranged in order to fill the gap below the left-most block in ${\rm Ap}_p(A)$. The other elements of ${\rm Ap}_{p-1}(A)$ are shifted directly to the right block by one in ${\rm Ap}_p(A)$.

In Table \ref{tb:g3system400}, $\ctext{$n$}$ denotes the area of elements in ${\rm Ap}_n(A)$. Here, each $m_j^{(n)}$, satisfying $m_j^{(n)}\equiv j\pmod{a}$ ($0\le j\le a-1$), can be expressed in at least $n+1$ ways but $m_j^{(n)}-a$ in at most $n$ ways. 
Each element of ${\rm Ap}_3(A)$ existing in the second block to the fourth block corresponds to each element having the same residue of ${\rm Ap}_2(A)$ existing in the block immediately to the left thereof in a form of shifting up one step.
The $J_k(v)$ elements of ${\rm Ap}_3(A)$ existing over two rows (or one row) at the bottom of the first block correspond to the $J_k(v)$ elements with the same residue of ${\rm Ap}_2(A)$ at the top of the third block. 
Therefore, since all the elements in ${\rm Ap}_2(A)$ form a complete remainder system, so is ${\rm Ap}_3(A)$.
It can be shown that all the elements of ${\rm Ap}_3(A)$ have exactly four ways of expressing in terms of $a$, $v a+b$, and $v a J_{k-1}(v)+b J_k(v)$. 
Within each region of ${\rm Ap}_3(A)$, one of the two leftmost (lower left) elements $t_{r-1,q+3}$ and $t_{J_k(v)-1,q+2}$ is the largest, so by comparing these elements, the largest element of ${\rm Ap}_3(A)$ can be determined. 

\begin{table}[htbp]
  \centering
\scalebox{0.7}{
\begin{tabular}{cccccccccccccccc}
\multicolumn{1}{|c}{}&&&\multicolumn{1}{c|}{}&&&&\multicolumn{1}{c|}{}&&&&\multicolumn{1}{c|}{}&&&&\multicolumn{1}{c|}{}\\ 
\cline{1-2}\cline{3-4}\cline{5-6}\cline{7-8}\cline{9-10}\cline{11-12}\cline{13-14}\cline{15-16}
\multicolumn{1}{|c}{}&&&\multicolumn{1}{c|}{}&&&&\multicolumn{1}{c|}{}&&&&\multicolumn{1}{c|}{}&&&&\multicolumn{1}{c|}{}\\ 
\multicolumn{1}{|c}{}&&&\multicolumn{1}{c|}{}&&&&\multicolumn{1}{c|}{}&&$\ctext{2}$&&\multicolumn{1}{c|}{}&&$\ctext{3}$&&\multicolumn{1}{c|}{}\\ 
\cline{14-15}\cline{15-16}
\multicolumn{1}{|c}{}&&&\multicolumn{1}{c|}{}&&$\ctext{1}$&&\multicolumn{1}{c|}{}&&&&\multicolumn{1}{c|}{}&\multicolumn{1}{c|}{}&&&\multicolumn{1}{c|}{}\\ 
\cline{10-11}\cline{12-13}
\multicolumn{1}{|c}{}&$\ctext{0}$&&\multicolumn{1}{c|}{}&&&&\multicolumn{1}{c|}{}&\multicolumn{1}{c|}{}&&$\ctext{3}$&\multicolumn{1}{c|}{}&&&&\multicolumn{1}{c|}{}\\ 
\cline{6-6}\cline{7-8}\cline{9-10}\cline{11-12}
\multicolumn{1}{|c}{}&&&\multicolumn{1}{c|}{}&\multicolumn{1}{c|}{}&&$\ctext{2}$&\multicolumn{1}{c|}{}&\multicolumn{1}{c|}{$\ctext{3}$}&&&\multicolumn{1}{c|}{}&&&&\multicolumn{1}{c|}{}\\ 
\cline{2-3}\cline{4-5}\cline{6-7}\cline{8-9}
\multicolumn{1}{|c|}{}&&$\ctext{1}$&\multicolumn{1}{c|}{}&\multicolumn{1}{c|}{$\ctext{2}$}&&$\ctext{3}$&\multicolumn{1}{c|}{}&&&&\multicolumn{1}{c|}{}&&&&\multicolumn{1}{c|}{}\\ 
\cline{1-2}\cline{3-4}\cline{5-6}\cline{7-8}
\multicolumn{1}{|c|}{$\ctext{1}$}&&$\ctext{2}$&\multicolumn{1}{c|}{}&\multicolumn{1}{c|}{$\ctext{3}$}&&&\multicolumn{1}{c|}{}&&&&\multicolumn{1}{c|}{}&&&&\multicolumn{1}{c|}{}\\ 
\cline{1-2}\cline{3-4}\cline{4-5}
\multicolumn{1}{|c|}{$\ctext{2}$}&&$\ctext{3}$&\multicolumn{1}{c|}{}&&&&\multicolumn{1}{c|}{}&&&&\multicolumn{1}{c|}{}&&&&\multicolumn{1}{c|}{}\\ 
\cline{1-2}\cline{3-4}
\multicolumn{1}{|c|}{$\ctext{3}$}&&&\multicolumn{1}{c|}{}&&&&\multicolumn{1}{c|}{}&&&&\multicolumn{1}{c|}{}&&&&\multicolumn{1}{c|}{}\\ 
\cline{1-1}
\multicolumn{1}{|c}{}&&&\multicolumn{1}{c|}{}&&&&\multicolumn{1}{c|}{}&&&&\multicolumn{1}{c|}{}&&&&\multicolumn{1}{c|}{}
\end{tabular}
} 
  \caption{${\rm Ap}_p(A)$ ($p=0,1,2,3$) for $q\ge 3$}
  \label{tb:g3system400}
\end{table}

Such a structure of ${\rm Ap}_p(A)$ continues as long as $p\le\fl{a/J_k(v)}=q$. Eventually, the largest element in ${\rm Ap}_p(A)$ is $t_{r-1,q+p}$ or $t_{J_k(v)-1,q+p-1}$. However, when $p=\fl{a/J_k(v)}+1$, this kind of regularity is broken. Therefore, regularity cannot be maintained even for the largest value of ${\rm Ap}_p(A)$ too. Therefore, Theorem \ref{th1} is proved. 
Table \ref{tb:gpsystemr=p} shows the arrangement of the $p$-Ap\'ery sets ($p=0,1,\dots,5$) when $\fl{a/J_k(v)}=5$. One can see that there will be a deficiency in the arrangement of the lower left for $6$-Ap\'ery set.

\begin{table}[htbp]
  \centering
\scalebox{0.7}{
\begin{tabular}{cccccccccccccccccccccccc}
\multicolumn{1}{|c}{}&&&\multicolumn{1}{c|}{}&&&&\multicolumn{1}{c|}{}&&&&\multicolumn{1}{c|}{}&&&&\multicolumn{1}{c|}{}&&&&\multicolumn{1}{c|}{}&&&&\multicolumn{1}{c|}{}\\ 
\cline{1-2}\cline{3-4}\cline{5-6}\cline{7-8}\cline{9-10}\cline{11-12}\cline{13-14}\cline{15-16}\cline{17-18}\cline{19-20}\cline{21-22}\cline{23-24}
\multicolumn{1}{|c}{}&&&\multicolumn{1}{c|}{}&&&&\multicolumn{1}{c|}{}&&&&\multicolumn{1}{c|}{}&&&&\multicolumn{1}{c|}{}&&$\ctext{4}$&&\multicolumn{1}{c|}{}&\multicolumn{1}{c|}{$\ctext{5}$}&&&\multicolumn{1}{c|}{}\\ 
\cline{18-18}\cline{19-20}\cline{21-21}
\multicolumn{1}{|c}{}&&&\multicolumn{1}{c|}{}&&&&\multicolumn{1}{c|}{}&&$\ctext{2}$&&\multicolumn{1}{c|}{}&&$\ctext{3}$&&\multicolumn{1}{c|}{}&\multicolumn{1}{c|}{}&&$\ctext{5}$&\multicolumn{1}{c|}{}&&&&\multicolumn{1}{c|}{}\\ 
\cline{14-15}\cline{15-16}\cline{17-18}\cline{19-20}
\multicolumn{1}{|c}{}&&&\multicolumn{1}{c|}{}&&$\ctext{1}$&&\multicolumn{1}{c|}{}&&&&\multicolumn{1}{c|}{}&\multicolumn{1}{c|}{}&&$\ctext{4}$&\multicolumn{1}{c|}{}&\multicolumn{1}{c|}{$\ctext{5}$}&&&\multicolumn{1}{c|}{}&&&&\multicolumn{1}{c|}{}\\ 
\cline{10-11}\cline{12-13}\cline{14-15}\cline{16-17}
\multicolumn{1}{|c}{}&$\ctext{0}$&&\multicolumn{1}{c|}{}&&&&\multicolumn{1}{c|}{}&\multicolumn{1}{c|}{}&&$\ctext{3}$&\multicolumn{1}{c|}{}&\multicolumn{1}{c|}{$\ctext{4}$}&&$\ctext{5}$&\multicolumn{1}{c|}{}&&&&\multicolumn{1}{c|}{}&&&&\multicolumn{1}{c|}{}\\ 
\cline{6-6}\cline{7-8}\cline{9-10}\cline{11-12}\cline{13-14}\cline{15-16}
\multicolumn{1}{|c}{}&&&\multicolumn{1}{c|}{}&\multicolumn{1}{c|}{}&&$\ctext{2}$&\multicolumn{1}{c|}{}&\multicolumn{1}{c|}{$\ctext{3}$}&&$\ctext{4}$&\multicolumn{1}{c|}{}&\multicolumn{1}{c|}{$\ctext{5}$}&&&\multicolumn{1}{c|}{}&&&&\multicolumn{1}{c|}{}&&&&\multicolumn{1}{c|}{}\\ 
\cline{2-3}\cline{4-5}\cline{6-7}\cline{8-9}\cline{10-11}\cline{12-13}
\multicolumn{1}{|c|}{}&&$\ctext{1}$&\multicolumn{1}{c|}{}&\multicolumn{1}{c|}{$\ctext{2}$}&&$\ctext{3}$&\multicolumn{1}{c|}{}&\multicolumn{1}{c|}{$\ctext{4}$}&&$\ctext{5}$&\multicolumn{1}{c|}{}&&&&\multicolumn{1}{c|}{}&&&&\multicolumn{1}{c|}{}&&&&\multicolumn{1}{c|}{}\\ 
\cline{1-2}\cline{3-4}\cline{5-6}\cline{7-8}\cline{9-10}\cline{11-12}
\multicolumn{1}{|c|}{$\ctext{1}$}&&$\ctext{2}$&\multicolumn{1}{c|}{}&\multicolumn{1}{c|}{$\ctext{3}$}&&$\ctext{4}$&\multicolumn{1}{c|}{}&\multicolumn{1}{c|}{$\ctext{5}$}&&&\multicolumn{1}{c|}{}&&&&\multicolumn{1}{c|}{}&&&&\multicolumn{1}{c|}{}&&&&\multicolumn{1}{c|}{}\\ 
\cline{1-2}\cline{3-4}\cline{5-6}\cline{7-8}\cline{8-9}
\multicolumn{1}{|c|}{$\ctext{2}$}&&$\ctext{3}$&\multicolumn{1}{c|}{}&\multicolumn{1}{c|}{$\ctext{4}$}&&$\ctext{5}$&\multicolumn{1}{c|}{}&&&&\multicolumn{1}{c|}{}&&&&\multicolumn{1}{c|}{}&&&&\multicolumn{1}{c|}{}&&&&\multicolumn{1}{c|}{}\\ 
\cline{1-2}\cline{3-4}\cline{5-6}\cline{7-8}
\multicolumn{1}{|c|}{$\ctext{3}$}&&$\ctext{4}$&\multicolumn{1}{c|}{}&\multicolumn{1}{c|}{$\ctext{5}$}&&&\multicolumn{1}{c|}{}&&&&\multicolumn{1}{c|}{}&&&&\multicolumn{1}{c|}{}&&&&\multicolumn{1}{c|}{}&&&&\multicolumn{1}{c|}{}\\ 
\cline{1-2}\cline{3-4}\cline{4-5}
\multicolumn{1}{|c|}{$\ctext{4}$}&&$\ctext{5}$&\multicolumn{1}{c|}{}&&&&\multicolumn{1}{c|}{}&&&&\multicolumn{1}{c|}{}&&&&\multicolumn{1}{c|}{}&&&&\multicolumn{1}{c|}{}&&&&\multicolumn{1}{c|}{}\\ 
\cline{1-2}\cline{3-4}
\multicolumn{1}{|c|}{$\ctext{5}$}&\multicolumn{1}{c}{\cellcolor[gray]{0.8}$\ctext{6}$}&&\multicolumn{1}{c|}{}&&&&\multicolumn{1}{c|}{}&&&&\multicolumn{1}{c|}{}&&&&\multicolumn{1}{c|}{}&&&&\multicolumn{1}{c|}{}&&&&\multicolumn{1}{c|}{}\\ 
\cline{1-1}
\multicolumn{1}{|c}{}&&&\multicolumn{1}{c|}{}&&&&\multicolumn{1}{c|}{}&&&&\multicolumn{1}{c|}{}&&&&\multicolumn{1}{c|}{}&&&&\multicolumn{1}{c|}{}&&&&\multicolumn{1}{c|}{}
\end{tabular}
} 
  \caption{${\rm Ap}_p(A)$ ($p=\fl{a/J_k(v)}$)}
  \label{tb:gpsystemr=p}
\end{table}

See \cite{KP,KY} etc. for a detailed explanation that the elements located within the entire specified areas actually constitute the elements of the $p$-Ap\'ery set. That is, they form a complete residue system modulo $a$, and each element is represented by $a$, $v a+b$, $v a J_{k-1}(v)+b J_k(v)$ in at least $p+1$ ways. The rough structure is similar to that in \cite{KP,KY}, though the structures of the $p$-Ap\'ery set in other cases are not necessarily similar or have not been known yet.

\section{$p$-genus}  

The elements of ${\rm Ap}_p(A)$ in the area of the $2 p$ staircase parts are 
\begin{align*} 
&t_{0,q+p},\dots,t_{r-1,q+p},\quad t_{r,q+p-1},\dots,t_{J_k(v)-1,q+p-1},\\
&t_{J_k(v),q+p-2},\dots,t_{J_k(v)+r-1,q+p-2},\quad t_{J_k(v)+r,q+p-3},\dots,t_{2 J_k(v)-1,q+p-3},\\
&t_{2 J_k(v),q+p-4},\dots,t_{2 J_k(v)+r-1,q+p-4},\quad t_{2 J_k(v)+r,q+p-5},\dots,t_{3 J_k(v)-1,q+p-5},\\
&\dots\\
&t_{(p-1)J_k(v),q-p+2},\dots,t_{(p-1)J_k(v)+r-1,q-p+2},\quad t_{(p-1)J_k(v)+r,q-p+1},\dots,t_{p J_k(v)-1,q-p+1}\\
\end{align*} 
in order from the lower left, and the elements of ${\rm Ap}_p(A)$ in the right-most main area are 
\begin{align*} 
&t_{p J_k(v),0},\dots,t_{p J_k(v)+r-1,0},\quad t_{p J_k(v)+r,0},\dots,t_{(p+1)J_k(v)-1,0},\\
&t_{p J_k(v),1},\dots,t_{p J_k(v)+r-1,1},\quad t_{p J_k(v)+r,1},\dots,t_{(p+1)J_k(v)-1,1},\\
&\dots\\
&t_{p J_k(v),q-p-1},\dots,t_{p J_k(v)+r-1,q-p-1},\quad t_{p J_k(v)+r,q-p-1},\dots,t_{(p+1)J_k(v)-1,q-p-1},\\
&t_{p J_k(v),q-p},\dots,t_{p J_k(v)+r-1,q-p}\,. 
\end{align*}   
Hence, by $a=q J_k(v)+r$, we have 
\begin{align}
&\sum_{w\in{\rm Ap}_p(A)}w\notag\\
&=\sum_{i=0}^{p-1}\sum_{j=0}^{r-1}t_{i J_k(v)+j,q+p-2 i}+\sum_{i=0}^{p-1}\sum_{j=0}^{J_k(v)-r-1}t_{i J_k(v)+r+j,q+p-2 i-1}\notag\\
&\quad +\sum_{i=0}^{J_k(v)-1}\sum_{j=0}^{q-p-1}t_{p J_k(v)+i,j}+\sum_{i=0}^{r-1}t_{p J_k(v)+i,q-p}\notag\\
&=\frac{a}{2}\bigl(-v(a-r^2)+b(a-1)+v\bigl((a+r)q-a+r\bigr)J_{k-1}(v)+v(a-r)J_k(v)\bigr)\notag\\
&\quad +\frac{p}{2}a J_k(v)\bigl(2(v a+b)-v(J_k(v)-J_{k-1}(v)\bigr)
-\frac{p^2}{2}a v J_k(v)\bigl(J_k(v)-J_{k-1}(v)\bigr)\,.
\label{eq:sum-w}
\end{align}
Thus, by Lemma \ref{lem-mp} (\ref{mp-n}), we obtain that  
\begin{align*}  
&n_p(a,v a+b,a J_{k-1}(v)+b J_k(v))\\
&=\frac{1}{a}\sum_{w\in{\rm Ap}_p(A)}w-\frac{a-1}{2}\\
&=\frac{1}{2}\bigl(-v(a-r^2)+b(a-1)+v\bigl((a+r)q-a+r\bigr)J_{k-1}(v)+v(a-r)J_k(v)\bigr)\\
&\quad +\frac{p}{2}a J_k(v)\bigl(2(v a+b)-v(J_k(v)-J_{k-1}(v)\bigr)\\
&\quad -\frac{p^2}{2}a v J_k(v)\bigl(J_k(v)-J_{k-1}(v)\bigr)
 -\frac{a-1}{2}\\
&=\frac{1}{2}\biggl(-v(a-r^2)+(a-1)(b-1)+\frac{v(a+r)(a-r)J_{k-1}(v)}{J_k(v)}\\
&\qquad +v(a-r)\bigl(J_k(v)-J_{k-1}(v)\bigr)\biggr)\\
&\quad +\frac{p}{2}a J_k(v)\bigl(2(v a+b)-v(J_k(v)-J_{k-1}(v)\bigr)
-\frac{p^2}{2}a v J_k(v)\bigl(J_k(v)-J_{k-1}(v)\bigr)\,. 
\end{align*}

\begin{theorem}
Let $a$ and $b$ be coprime integers. Then, 
for a positive integer $k\ge 3$, and $0\le p\le\fl{a/J_k(v)}$ we have 
\begin{align*}
&n_p\bigl(a,a+b,a J_{k-1}(v)+b J_k(v)\bigr)\\
&=\frac{1}{2}\biggl(-v(a-r^2)+(a-1)(b-1)+\frac{v(a+r)(a-r)J_{k-1}(v)}{J_k(v)}\\
&\qquad +v(a-r)\bigl(J_k(v)-J_{k-1}(v)\bigr)\biggr)\\
&\quad +\frac{p}{2}a J_k(v)\bigl(2(v a+b)-v(J_k(v)-J_{k-1}(v)\bigr)
-\frac{p^2}{2}a v J_k(v)\bigl(J_k(v)-J_{k-1}(v)\bigr)\,, 
\end{align*} 
where $r=a-\fl{a/J_k(v)}J_k(v)$. 
\label{th2}
\end{theorem}

\subsection{$p$-Sylvester sum}  

In this subsection, we shall show a closed formula for the Sylvester sum.   
By $a=q J_k(v)+r$, we have 
\begin{align*}
&\sum_{w\in{\rm Ap}_p(A)}w^2\\
&=\sum_{i=0}^{p-1}\sum_{j=0}^{r-1}(t_{i J_k(v)+j,q+p-2 i})^2+\sum_{i=0}^{p-1}\sum_{j=0}^{J_k(v)-r-1}(t_{i J_k(v)+r+j,q+p-2 i-1})^2\\
&\quad +\sum_{i=0}^{J_k(v)-1}\sum_{j=0}^{q-p-1}(t_{p J_k(v)+i,j})^2+\sum_{i=0}^{r-1}(t_{p J_k(v)+i,q-p})^2\\
&=\frac{a}{6}\biggl(
(v a+b)^2+(v a+b)(2 v r^2-3 a b-3 v r^2)+2 a b(a b+3 v r^2)\\
&\qquad +v^2 a\bigl(2 q^2(a+2 r)-(a-r)(3 q-1)\bigr)J_{k-1}(v)^2\\
&\qquad -3 v(a-r)(v a+b-b(a-r)\bigr)J_k(v)+v(a-r)(2 v a+b)J_k(v)^2\\
&\qquad 
+v\bigl(3(a-r)((a-r)(v a-b)+v a+b)+q(a+r)(4 a b-3(v a+b))\\
&\qquad +2 q r^2(3 v a-b)-(a-r)(3 v a+b)J_k(v)\bigr)J_{k-1}(v)
\biggr)\\ 
&\quad +\frac{a p}{6}\biggl(
6\bigl((v a+b)(a b+v r^2)-(v a+b)^2-b^2\bigr)J_k(v)\\
&\qquad +3 v\bigl((v a+b)(2 r+1)+2 v a^2\bigr)J_k(v)^2-v(2 v a-b)J_k(v)^3\\
&\qquad +v\bigl(6(a^2-r^2)(a v+b)+3((v a+b)(2 r-1)-4 v a^2)J_k(v)\\
&\qquad +(5 v a-b)J_k(v)^2\bigr)J_{k-1}(v)+3 v^2 a\bigl(2 a-J_k(v)\bigr)J_{k-1}(v)^2
\biggr)\\ 
&\quad+\frac{a p^2}{2}\biggl(\bigl(v(v a+b)(2 a+1)+2 b^2-v(2 v a+b)J_k(v)\bigr)J_k(v)^2\\
&\qquad -v\bigl(2 a(v a-b)+a v+b\bigr)J_{k-1}(v)J_k(v)+v^2 a J_{k-1}(v)^2(2 a-J_k(v)^2)\biggr)\\ 
&\quad -\frac{2 a v(v a+b)p^3}{3}\bigl(J_{k}(v)-J_{k-1}(v)\bigr)(J_k(v))^2\,. 
\end{align*} 
Thus, by Lemma \ref{lem-mp} (\ref{mp-s}),  
together with $\sum_{w\in{\rm Ap}_p(A)}w$ in (\ref{eq:sum-w}), we obtain that 
\begin{align*}  
&s_p\bigl(a,v a+b,v a J_{k-1}(v)+b J_k(v)\bigr)\\
&=\frac{1}{2 a}\sum_{w\in{\rm Ap}_p(A)}w^2-\frac{1}{2}\sum_{w\in{\rm Ap}_p(A)}w+\frac{a^2-1}{12}\\
&=\frac{1}{12}\biggl(
(a-r)v(2 a v+b)J_k(v)^2-(a-r)v\bigl(3(v a+b-b(a-r)+a)\\
&\qquad +(3 v a+b)J_{k-1}(v)\bigr)J_k(v)+v^2 a(a-r)J_{k-1}(v)^2
\\
&\qquad +3(a-r)v\bigl(v a+b+(a-r)(v a-b)+a\bigr)J_{k-1}(v)\\
&\qquad +3\bigl(v r^2(v a+b)+v a r^2(2 b-1)-v a^2(b-1)-a b(a+b-1)\bigr)\\
&\qquad +(v a+b)^2+(a^2-1)+2\bigl(v r^3(v a-b)+a^2 b^2\bigr)\\ 
&\qquad +\frac{2 v^2 a(a-r)(a^2+a r-2 r^2)J_{k-1}(v)^2}{J_k(v)^2}-\frac{3 v^2 a(a-r)^2 J_{k-1}(v)^2}{J_k(v)}\\
&\qquad +\frac{v(a-r)J_{k-1}(v)}{J_k(v)}\bigl(4 a b(2 a+3 r)+2 r^2(3 v-b)\\
&\qquad\quad -3(a^2+r(v+b+1))-3 a(v a+b)\bigr)
\biggr)\\ 
&\quad +\frac{p}{12}\biggl(-v(2 v a-b)J_k(v)^3+3 v\bigl(-(2 r-1)(v a+b)+a(2 v a+1)\bigr)J_k(v)^2\\
&\qquad +6(v a+b)\bigl(v r^2-(v a+b)+a(b-1)\bigr)J_k(v)\\
&\qquad +v(5 v a-b)J_{k-1}(v)J_k(v)^2\\
&\qquad +3 v\bigl(-(2 r-1)(v a+b)+a(4 v a+1)\bigr)J_{k-1}(v)J_k(v)\\
&\qquad +6 v(a^2-r^2)(v a+b)J_{k-1}(v)+3 v^2 a\bigl(2 a-J_k(v)\bigr)J_{k-1}(v)^2
\biggr)\\ 
&\quad+\frac{p^2}{4}\biggl(v\bigl(a(2 v a-2 b+1)+v a+b-(3 v a+b)J_k(v)\bigr)J_{k-1}(v)J_k(v)\\
&\qquad +\bigl(v(2 a+1)(v a+b)+v a+2 b^2-v(2 v a+b)J_k(v)\bigr)J_k(v)^2\\
&\qquad +v^2 a J_{k-1}(v)^2 J_k(v)^2\biggr)\\ 
&-\frac{v(v a+b)p^3}{3}\bigl(J_{k}(v)-J_{k-1}(v)\bigr)J_k(v)^2\,. 
\end{align*} 
Here, again $q=\fl{a/J_k(v)}$ and $r=a-q J_k(v)$.

\section{Jacobsthal-Lucas polynomials}  

The same discussion as Jacobsthal polynomials can be applied to Jacobsthal-Lucas polynomials $j_n(v)$ as it is.  Here, $j_n(v)=j_{n-1}(v)+j_{n-2}(v)$ ($n\ge 2$) with $j_0(v)=2$ and $j_1(v)=1$ (see, e.g., \cite[Chapter 44]{Koshy2}). When $v=1$, $L_n=j_n(1)$ are Lucas numbers. When $v=2$, $j_n=j_n(2)$ are Jacobsthal-Lucas numbers. 
Similarly, determine integers $q$ and $r$ by $a=q j_k(v)+r$ with $0\le r<j_k(v)$. If $a=j_i(v)$ and $b=j_{i+1}(v)$, then the numerical semigroup $\ang{j_i(v),j_{i+2}(v),j_{i+k}(v)}$ in \cite{KP} can be reduced as a special case.  

\begin{theorem}
Let $a$ and $b$ be positive integers with $\gcd(a,b)=1$ and $a\ge 3$. Then, 
for a positive integer $k\ge 3$, and $0\le p\le\fl{a/j_k(v)}$ we have  
\begin{align*}
&g_p\bigl(a,v a+b,v a j_{k-1}(v)+b j_k(v)\bigr)\\
&=\begin{cases}
(a-1)b+a\bigl(v(r-1)-1\bigr)+\frac{v a(a-r)j_{k-1}(v)}{j_k(v)}+p\bigl(v a j_{k-1}(v)+b j_k(v)\bigr)&\\
\qquad\qquad\qquad\text{if $a<j_k(v)$ or $(v a+b)r>v a\bigl(j_k(v)-j_{k-1}(v)\bigr)$};&\\ 
(a-r-1)b+v a(j_k(v)-j_{k-1}(v)-1)-a+\frac{v a(a-r)j_{k-1}(v)}{j_k(v)}&\\
+p\bigl(v a j_{k-1}(v)+b j_k(v)\bigr)&\\
\qquad\qquad\qquad\text{if $a\ge j_k(v)$ and $(v a+b)r<v a\bigl(j_k(v)-j_{k-1}(v)\bigr)$}\,,& 
\end{cases}
\end{align*} 
where $r=a-\fl{a/j_k(v)}j_k(v)$. 
\label{th1:lucas}
\end{theorem}  

\begin{theorem}
Let $a$ and $b$ be coprime integers. Then, 
for a positive integer $k\ge 3$, and $0\le p\le\fl{a/j_k(v)}$ we have 
\begin{align*}
&n_p\bigl(a,a+b,a j_{k-1}(v)+b j_k(v)\bigr)\\
&=\frac{1}{2}\biggl(-v(a-r^2)+(a-1)(b-1)+\frac{v(a+r)(a-r)j_{k-1}(v)}{j_k(v)}\\
&\qquad +v(a-r)\bigl(j_k(v)-j_{k-1}(v)\bigr)\biggr)\\
&\quad +\frac{p}{2}a j_k(v)\bigl(2(v a+b)-v(j_k(v)-j_{k-1}(v)\bigr)
-\frac{p^2}{2}a v j_k(v)\bigl(j_k(v)-j_{k-1}(v)\bigr)\,, 
\end{align*} 
where $r=a-\fl{a/j_k(v)}j_k(v)$. 
\label{th2:lucas}
\end{theorem}

\section{$p$-Hilbert series}   

There are some applications, due to the $p$-Ap\'ery set. One of them is on the {\it $p$-Hilbert series} (\cite{Ko-p}) of $S_p(A)$, which is defined by 
$$
H_p(A;x):=H(S_p;x)=\sum_{s\in S_p(A)}x^s\,. 
$$  
When $p=0$, the $0$-Hilbert series is the original Hilbert series, which plays the important role in numerical semigroup (see, e.g., \cite{ADG}).  
In addition, the $p$-gaps generating function is defined by   
$$
\Psi_p(A;x)=\sum_{s\in\mathbb N_0\backslash S_p(A)}x^s\,,  
$$  
satisfying $H_p(A;x)+\Psi_p(A;x)=1/(1-x)$ ($|x|<1$). 
Moreover, the same arguments of Chapter 5 in \cite{ADG}, we can express the $p$-Hilbert series as  
\begin{equation} 
H_p(A;x)=\frac{1}{1-x^a}\sum_{w\in {\rm Ap}_p(A;a)}x^w\,, 
\label{eq:p-hilbert}
\end{equation} 
where $a=\min\{A\}$.     

When $A=\{a,v a+b,v a J_{k-1}(v)+b J_k(v)\}$, similarly to (\ref{eq:sum-w}), we have  
\begin{align*}
&\sum_{w\in{\rm Ap}_p(A)}x^w\\
&=\sum_{i=0}^{p-1}\sum_{j=0}^{r-1}x^{t_{i J_k(v)+j,q+p-2 i}}+\sum_{i=0}^{p-1}\sum_{j=0}^{J_k(v)-r-1}x^{t_{i J_k(v)+r+j,q+p-2 i-1}}\\
&\quad +\sum_{i=0}^{J_k(v)-1}\sum_{j=0}^{q-p-1}x^{t_{p J_k(v)+i,j}}+\sum_{i=0}^{r-1}x^{t_{p J_k(v)+i,q-p}}\\
&=\frac{(1-x^{r(v a+b)})(x^{2 p(v a J_{k-1}(v)+b J_k(v))}-x^{p(v a+b)J_k(v)})}{x^{(p-q-2)(v a J_{k-1}(v)+b J_k(v))}(1-x^{v a+b})(x^{2(v a J_{k-1}(v)+b J_k(v))}-x^{(v a+b)J_k(v)})}\\
&\quad +\frac{(x^{r(v a+b)}-x^{(v a+b)J_k(v)})(x^{2 p(v a J_{k-1}(v)+b J_k(v))}-x^{p(v a+b)J_k(v)})}{x^{(p-q-1)(v a J_{k-1}(v)+b J_k(v))}(1-x^{v a+b})(x^{2(v a J_{k-1}(v)+b J_k(v))}-x^{(v a+b)J_k(v)})}\\
&\quad +\frac{(x^{p(v a+b)J_k(v)}-x^{(p+1)(v a+b)J_k(v)})(1-x^{(q-p)(v a J_{k-1}(v)+b J_k(v)))})}{(1-x^{v a+b})(1-x^{v a J_{k-1}(v)+b J_k(v)})}\\ 
&\quad +\frac{x^{p v a(J_k(v)-J_{k-1}(v))+q(v a J_{k-1}+b J_k(v))}(1-x^{r(v a+b)})}{1-x^{v a+b}}\,. 
\end{align*} 
Therefore, by (\ref{eq:p-hilbert})  
\begin{align*}
&H_p\bigl(a,v a+b,v a J_{k-1}(v)+b J_k(v);x\bigr)\\
&=\frac{1}{1-x^a}\biggl(\frac{(1-x^{r(v a+b)})(x^{2 p(v a J_{k-1}(v)+b J_k(v))}-x^{p(v a+b)J_k(v)})}{x^{(p-q-2)(v a J_{k-1}(v)+b J_k(v))}(1-x^{v a+b})(x^{2(v a J_{k-1}(v)+b J_k(v))}-x^{(v a+b)J_k(v)})}\\
&\quad +\frac{(x^{r(v a+b)}-x^{(v a+b)J_k(v)})(x^{2 p(v a J_{k-1}(v)+b J_k(v))}-x^{p(v a+b)J_k(v)})}{x^{(p-q-1)(v a J_{k-1}(v)+b J_k(v))}(1-x^{v a+b})(x^{2(v a J_{k-1}(v)+b J_k(v))}-x^{(v a+b)J_k(v)})}\\
&\quad +\frac{(x^{p(v a+b)J_k(v)}-x^{(p+1)(v a+b)J_k(v)})(1-x^{(q-p)(v a J_{k-1}(v)+b J_k(v)))})}{(1-x^{v a+b})(1-x^{v a J_{k-1}(v)+b J_k(v)})}\\ 
&\quad +\frac{x^{p v a(J_k(v)-J_{k-1}(v))+q(v a J_{k-1}+b J_k(v))}(1-x^{r(v a+b)})}{1-x^{v a+b}}\biggr)\,. 
\end{align*}

\section*{Acknowledgement}  

The most of this paper has been completed when the first named author stayed in Delhi in February 2023. He thanks the discussions and hospitality by the second and third named authors.


\begin{thebibliography}{99} 

\bibitem{Apery}  
R. Ap\'ery,  {\em 
Sur les branches superlin\'eaires des courbes alg\'ebriques}, 
C. R. Acad. Sci. Paris {\bf 222} (1946), 1198--1200. 

\bibitem{ADG}  
A. Assi, M. D'Anna and P. A. Garcia-Sanchez, {\em 
Numerical semigroups and applications}, 
2nd extended and revised edition, 
RSME Springer Series 3. Cham: Springer (2020).  

\bibitem{SKT}   
S. S. Batra, N. Kumar and A. Tripathi, {\em  
On a linear Diophantine problem involving the Fibonacci and Lucas sequences}, 
Integers {\bf 15} (2015), Paper No. A26, 12 pp. 

\bibitem{Binner20}
D. S. Binner, {\em  
The number of solutions to $ax+by+cz=n$ and its relation to quadratic residues}, J. Integer Seq. {\bf 23} (2020), No. 6, Art. 20.6.5, 19 pp. 

\bibitem{bs62} 
A. Brauer and B. M. Shockley, {\em  
On a problem of Frobenius}, 
J. Reine. Angew. Math. {\bf 211} (1962), 215--220.  

\bibitem{BDFHKMRSS} 
A. Brown, E. Dannenberg, J. Fox, J. Hanna, K. Keck, A. Moore, Z. Robbins, B. Samples and J. Stankewicz, {\em 
On a generalization of the Frobenius number}, 
arXiv:1001.0207 (2010). 

\bibitem{bs93} 
T. C. Brown and P. J. Shiue,  {\em  
A remark related to the Frobenius problem},  
Fibonacci Quart. {\bf 31} (1993), 32--36.  

\bibitem{Cayley}  
A. Cayley, {\em On a problem of double partitions}, 
Philos. Mag. {\bf XX} (1860), 337--341. 

\bibitem{cu90} 
F. Curtis, {\em 
On formulas for the Frobenius number of a numerical semigroup}, 
Math. Scand. {\bf 67} (1990), 190--192. 

\bibitem{Fel06}  
L. G. Fel, {\em  
Frobenius problem for semigroups $S(d_1,d_2,d_3)$}, 
Funct. Anal. Other Math. {\bf 1} (2006), no. 2, 119--157.

\bibitem{FS11} 
L. Fukshansky and A. Schurmann, {\em 
Bounds on generalized Frobenius numbers}, 
Eur. J. Comb. {\bf 32} (2011), No. 3, 361--368.  

\bibitem{Johnson60}
S. M. Johnson, {\em  
A linear diophantine problem}, 
Canadian J. Math. {\bf 12} (1960), 390--398. 

\bibitem{Ko03}  
T. Komatsu, {\em  
On the number of solutions of the Diophantine equation of Frobenius--General case}, 
Math. Commun. {\bf 8} (2003), 195--206. 

\bibitem{Ko22} 
T. Komatsu, {\em 
Sylvester power and weighted sums on the Frobenius set in arithmetic progression}, 
Discrete Appl. Math. {\bf 315} (2022), 110--126. 

\bibitem{Ko22a} 
T. Komatsu, {\em 
The Frobenius number for sequences of triangular numbers associated with number of solutions},  
Ann. Comb. {\bf 26} (2022) 757--779. 

\bibitem{Ko22b}
T. Komatsu, {\em 
The Frobenius number associated with the number of representations for sequences of repunits}, 
C. R. Math., Acad. Sci. Paris {\bf 361} (2023), 73--89.   
https://doi.org/10.5802/crmath.394  

\bibitem{Ko-p} 
T. Komatsu, {\em 
On $p$-Frobenius and related numbers due to $p$-Ap\'ery set},  
arXiv:2111.11021v3 (2023).  

\bibitem{KP}  
T. Komatsu and C. Pita-Ruiz, {\em  
The Frobenius number for Jacobsthal triples associated with number of solutions}, 
Axioms {\bf 12} (2023), Article 98, 18 pp. 
https://doi.org/10.3390/axioms12020098 

\bibitem{KY}  
T. Komatsu and H. Ying, {\em  
The $p$-Frobenius and $p$-Sylvester numbers for Fibonacci and Lucas triplets},  
Math. Biosci. Eng. {\bf 20} (2023), No.2, 3455--3481. doi:10.3934/mbe.2023162

\bibitem{KZ0}  
T. Komatsu and Y. Zhang, {\em 
Weighted Sylvester sums on the Frobenius set}, 
Irish Math. Soc. Bull. {\bf 87} (2021), 21--29.   

\bibitem{KZ}  
T. Komatsu and Y. Zhang, {\em 
Weighted Sylvester sums on the Frobenius set in more variables}, 
Kyushu J. Math. {\bf 76} (2022), 163--175.   

\bibitem{Koshy2} 
Koshy, T. {\em  
Fibonacci and Lucas Numbers with Applications};  
Pure and Applied Mathematics (Hoboken); John Wiley \& Sons, Inc.: Hoboken, NJ, USA, 2019; Volume 2.

\bibitem{Matt09}  
G. L. Matthews, {\em 
Frobenius numbers of generalized Fibonacci semigroups},  
Landman, Bruce (ed.) et al., Combinatorial number theory. Proceedings of the 3rd `Integers Conference 2007', Carrollton, GA, USA, October 24--27, 2007. Berlin: Walter de Gruyter (ISBN 978-3-11-020221-2/hbk), 
Integers {\bf 9} Suppl., Article A9, 117--124 (2009).   

\bibitem{MAR} 
J. M. Marin, J. L. Ramirez Alfonsin and M. P. Revuelta,  {\em 
On the Frobenius number of Fibonacci numerical semigroups},  
Integers {\bf 7} (2007), A14, 7 pp. 

\bibitem{PT}  
P. Punyani and A. Tripathi, {\em 
On changes in the Frobenius and Sylvester numbers}, 
Integers {\bf 18B} (2018), Paper No. A8, 12 pp. 

\bibitem{RA05} 
J. L. Ramirez Alfonsin, {\em  
The Diophantine Frobenius problem}, 
Oxford Lecture Series in Mathematics and its Applications, 30. Oxford University Press, Oxford, 2005. 

\bibitem{RG1} 
J. C. Rosales and P. A. Garcia-Sanchez, {\em   
Finitely generated commutative monoids}, 
Nova Science Publishers, Inc., Commack, NY, 1999.   

\bibitem{RG2}  
J. C. Rosales and P. A. Garcia-Sanchez, {\em  
Numerical semigroups}, 
Developments in Mathematics, {\bf 20}. Springer, New York, 2009. 

\bibitem{RR18}  
A. M. Robles-P\'erez and J. C. Rosales,  {\em  
The Frobenius number for sequences of triangular and tetrahedral numbers}, 
J. Number Theory {\bf 186} (2018), 473--492.  

\bibitem{ro94} 
\O. J. R\o dseth, {\em 
A note on Brown and Shiue's paper on a remark related to the Frobenius problem},
Fibonacci Quart. {\bf 32}  (1994), 407--408. 

\bibitem{RGS04}
J. C. Rosales and P. A. Garcia-Sanchez, {\em 
Numerical semigroups with embedding dimension three}, 
Arch. Math. (Basel) {\bf 83} (2004), no. 6, 488--496. 

\bibitem{RBT2015}  
J. C. Rosales, M. B. Branco and D. Torr\~ao, {\em 
The Frobenius problem for Thabit numerical semigroups},  
J. Number Theory {\bf 155} (2015), 85--99. 

\bibitem{RBT2017}  
J. C. Rosales, M. B. Branco and D. Torr\~ao, {\em 
The Frobenius problem for Mersenne numerical semigroups},  
Math. Z. {\bf 286} (2017), 741--749.  

\bibitem{se77}  
E. S. Selmer, {\em  
On the linear diophantine problem of Frobenius},  
J. Reine Angew. Math. {\bf 293/294} (1977), 1--17.  

\bibitem{sy1857} 
J. J. Sylvester, {\em 
On the partition of numbers}, 
Quart. J. Pure Appl. Math. {\bf 1} (1857), 141--152.

\bibitem{sy1882}  
J. J. Sylvester, {\em 
On subinvariants, i.e. semi-invariants to binary quantics of an unlimited order},
Am. J. Math. {\bf 5} (1882), 119--136.

\bibitem{sy1884}  
J. J. Sylvester,  {\em  
Mathematical questions with their solutions},  
Educational Times {\bf 41} (1884), 21.  

\bibitem{tr00}  
A. Tripathi, {\em 
The number of solutions to $a x+b y=n$}, 
Fibonacci Quart. {\bf 38} (2000), 290--293.

\bibitem{tr08}  
A. Tripathi, {\em 
On sums of positive integers that are not of the form $a x+b y$}, 
Amer. Math. Monthly {\bf 115} (2008), 363--364.  

\bibitem{tr17} 
A. Tripathi, {\em 
Formulae for the Frobenius number in three variables}, 
J. Number Theory {\bf 170} (2017), 368--389.

\end{thebibliography}
\end{document}